\documentclass[reprint,
amsmath,amssymb,
aps,superscriptaddress,
]{revtex4-1}

\usepackage{amsmath,amsthm,amssymb}
\usepackage{graphicx}
\usepackage{color}
\usepackage{natbib}
\usepackage{subfigure}
\usepackage{multirow}

\def \deta {\pmb{\eta}}

\begin{document} 

\title{Analyzing synchronized clusters in neuron networks}

\author{Matteo Lodi}
\affiliation{DITEN, University of Genoa, Via Opera Pia 11a, I-16145, Genova, Italy}
\author{Fabio Della Rossa}
\affiliation{Mechanical Engineering Department, University of New Mexico, Albuquerque, NM 87131, USA}
\affiliation{DEIB, Politecnico di Milano, I-20133 Milan, Italy}
\author{Francesco Sorrentino}
\affiliation{Mechanical Engineering Department, University of New Mexico, Albuquerque, NM 87131, USA}
\author{Marco Storace}
\affiliation{DITEN, University of Genoa, Via Opera Pia 11a, I-16145, Genova, Italy}
\affiliation{Corresponding author: marco.storace@unige.it.}

\date{\today}
	
\begin{abstract}
The presence of synchronized clusters in neuron networks is a hallmark of information transmission and processing. The methods commonly used to study cluster synchronization in networks of coupled oscillators ground on simplifying assumptions, which often neglect key biological features of neuron networks. Here we propose a general framework to study presence and stability of synchronous clusters in more realistic models of neuron networks, characterized by the presence of delays, different kinds of neurons and synapses. 
Application of this framework to the directed network of the macaque cerebral cortex provides an interpretation key to explain known functional mechanisms emerging from the combination of anatomy and neuron dynamics.
The cluster synchronization analysis is carried out also by changing parameters and studying bifurcations. 
Despite some simplifications with respect to the real network, the obtained results are in good agreement with previously reported biological data.
\end{abstract}

\maketitle

Understanding the functional mechanisms of a given system/phenomenon and describing it through mathematical equations as simple as possible (according to the \textit{Occam's razor} principle) is the Holy Grail of modeling.
Among the others, neuron networks are the object of many studies due to their complex behaviors; understanding the functional mechanisms of information transmission and processing in these networks is one of the most difficult and fascinating challenges faced by the scientific community, at the crossroad between many disciplines.
The level of abstraction used to describe neuron networks can significantly change according to the modeling goals, complexity of the network to be modeled and background knowledge \cite{bassett:2018}. Consequently, the basic elements of the nervous system (neurons and synapses) are modeled by trading off accuracy and complexity \cite{herz:2006}. Neurons in the same network can be of different kinds and their synaptic connections, also of different kinds, can be either electrical or chemical, either excitatory or inhibitory, either directed or undirected, and may transmit signals with different delays. In this letter, we focus on deterministic models of these networks.

A commonly observed phenomenon in networks of neurons is the formation of \textit{synchronous clusters}, i.e., groups of neurons that fulfill some synchrony conditions \cite{kreiter:1996,maldonado:2000,glennon:2016}, usually expressed in terms of temporal correlation between neural signals. These clusters are strongly related to information transmission and processing \cite{bullmore:2009}. Recent efforts have been devoted to apply nonlinear dynamics concepts and network theory to the neuroscience context \cite{guevara:2017,bassett:2018}. However, deterministic models that study the presence and the stability of synchronized clusters in networks are based on simplifying assumptions such as identical neurons/synapses, weak interactions, absence of delays, undirected/diffusive connections. In this letter we propose a general method to analyze cluster synchronization (CS) in neuron networks with directed connections, delays, couplings that depend on both the presynaptic and the postsynaptic neurons, and different kinds of nodes and synapses. 
These networks can be described by the following set of dynamical equations, describing a multi-layer network \cite{boccaletti:2014},
($i = 1,\ldots,N$)
\begin{equation}
\dot x_i =  \tilde{f}_i(x_i(t)) + \sum_{k=1}^{L}\sigma^k \sum_{j=1}^{N} A_{ij}^k h^k(x_i(t),x_j(t-\delta_k)),
\label{eq:dyneq}
\end{equation}
where $x_i \in \mathbb{R}^n$ is the $n$-dimensional state vector of the $i$-th neuron, $\tilde{f}_i : \mathbb{R}^n \rightarrow \mathbb{R}^n$ is the vector field of the isolated $i$-th neuron, $\sigma^k \in \mathbb{R}$ is the coupling strength of the $k$-th kind of link, $A^k$ is the possibly weighted and directed coupling matrix (or adjacency matrix) that describes the connectivity of the network with respect to the $k$-th kind of link, for which the interaction between two generic cells $i$ and $j$ is described by the nonlinear function $h^k(\cdot) : \mathbb{R}^n \times \mathbb{R}^n \rightarrow \mathbb{R}^n$, and $\delta_k$ is the axon transmission delay characteristic of the $k$-th kind of link. 
For example, electrical synapses (gap junctions) are almost instantaneous, whereas the delay associated with transmission of a signal through a chemical synapse may be considerably longer.

A neuron model is described by a state vector $x_i$, whose first component $V_i$ typically represents the membrane potential of the neuron.
A synapse model can either neglect or include the neurotransmitter dynamics, therefore we can have instantaneous or dynamical synapses, respectively. In both cases, we assume that the synaptic coupling influences only the dynamics of $V_i$ and not of the other state variables contained in $x_i$: therefore, the first component of the vector $h^k(\cdot)$ is a scalar function (called \textit{activation function}) $a^k(V_i(t),x_j(t-\delta_k))$ and the remaining components are null.
For instantaneous synapses, the activation depends on the membrane potential of the pre- and post-synaptic neurons, therefore it can be expressed as $a^k(V_i(t),V_j(t-\delta_k))$. By contrast, for dynamical synapses the activation $a^k$ is a function of a state variable $s_j^k$ (in addition to $V_i$), whose dynamics usually depends on the pre-synaptic membrane potential $V_j$ (see Sec. 1 in \footnote{See Supplemental Material for further examples, mathematical details on the proposed method, and additional datasets.}). For this reason, all dynamical synapses of kind $k$ connecting the neuron $j$ with other neurons share the same state $s_j^k$, which can be added to vector $x_j$.

We further assume each individual node can be of one out of $M$ different types (with $M \leq N$): $\tilde{f}_i(x) = \tilde{f}_j(x)$ if $i$ and $j$ are of the same type, $\tilde{f}_i(x) \neq \tilde{f}_j(x)$ otherwise. Often, the difference (physical or functional) between two types of neurons is accounted for through a different value of one or more model parameters.
Within this general framework, where all oscillators can be different, if $M << N$ the vector fields $\tilde{f}_i$ are not all different, but belong to a restricted set of $M$ models.
Assuming that all node states share the same dimension $n$ is not restrictive: in the case of state vectors $x_i$ with different dimensions $n_i$, it is sufficient to define $\displaystyle{n = \max_{i} n_i}$ and set to 0 the components in excess [9].

Different from most models introduced in the literature, the set of equations (1) accounts for the following  realistic properties of neuron networks: (i) each synapse depends (algebraically in the case of instantaneous/fast synapses or dynamically in the case of slower synapses) on the state of both the pre-synaptic and the post-synaptic neuron, (ii) each synapse between two neurons is in general a direct connection that can be of different kinds (such as either chemical inhibitory/excitatory or electrical excitatory), and (iii) the transmission of information along synapses can be non-instantaneous, which may be due in part to local synaptic filtering of exchanged spikes, and in part to the distribution of the axonal transmission delays \cite{mattia:2019}. We wish to emphasize that current methods developed to analyze CS in complex networks \cite{pecora:2014,dellarossa:2020} are unable to handle features (i), (ii) and (iii) above.




\begin{figure*}[t!]
	\centering
	\includegraphics[width=\textwidth]{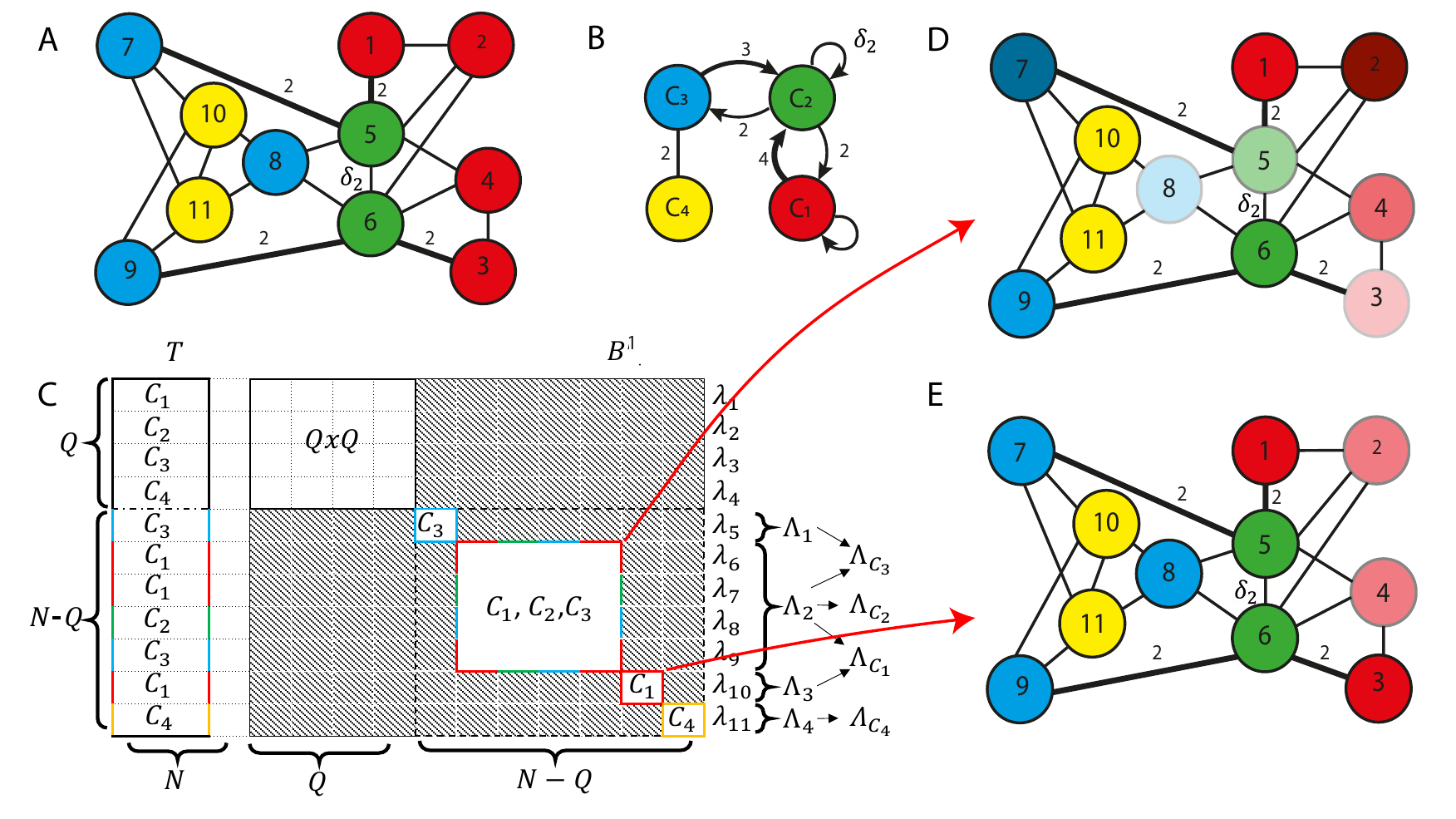}
	\caption{Example: (A) a network with $N=11$ nodes, $L=2$ kind of connections, and $Q=4$ clusters ($C_1 = \{1,2,3,4\},C_2 = \{5,6\},C_3 = \{7,8,9\},C_4 = \{10,11\}$) and (B) the structure of the corresponding matrices $T$ and $B^1$, illustrating their relation with the clusters. All connections are bi-directional and with weight 1, with the exception of the thick connection between nodes 5 and 7, which has weight 2. Network coloring (with a larger number of clusters) after the breaking of the red cluster if its loss of stability is due to the MLEs corresponding to either  (C) the multi-color sub-block or (D) the red sub-block of matrix $B^1$.}\label{fig:netwTmatrix}
\end{figure*}

Cluster synchronization of the system in Eq. \eqref{eq:dyneq} is defined as $x_i(t) = x_j(t)$ for any $t$ and for $i,j$ belonging to the same cluster of a certain partition. The set of the network nodes can be partitioned into equitable clusters (ECs), whose presence is necessary to achieve CS. Indeed, nodes in the same EC receive the same amount of weighted inputs of a certain type from the other ECs or from the EC itself.
The method we propose for the analysis of CS in networks modeled by Eq. \eqref{eq:dyneq} consists of three main steps: (S1) a coloring algorithm to find the $Q$ ECs $C_q$ ($q = 1,\ldots,Q$) of the network, corresponding to a clustering $\mathcal{C} = \{C_1,\ldots,C_Q\}$ (see the example network in Fig. \ref{fig:netwTmatrix}A, where $N=11$ and $Q=4$); (S2) a simplified dynamical model (called \textit{quotient network}) whose $Q$ nodes correspond to each one of the ECs (see Fig. \ref{fig:netwTmatrix}B, which is the quotient network corresponding to Fig. \ref{fig:netwTmatrix}A); (S3) an analysis of the cluster stability by linearizing Eq. \eqref{eq:dyneq} about a state corresponding to exact synchronization among all the nodes within each cluster. A detailed description of steps S1, S2, and S3 is provided in [9]. The main novelty of this method is the analysis S3, which is tailored to Eq. \eqref{eq:dyneq} following, \textit{mutatis mutandis}, the guidelines defined in \cite{pecora:2014,dellarossa:2020}. A key step of this analysis is the construction of the (unique) matrix $T$ that transforms the coupling matrices $A^k$ into block diagonal matrices, $B^k=T A^k T^T$. This corresponds to a change of perturbation coordinates that converts the node coordinate system to the \textit{irreducible representation} \cite{golubitsky:2012,pecora:2014,dellarossa:2020} coordinate system, thus evidencing the interdependencies among the perturbation components. For undirected networks, the $N \times N$ matrix $T$ can be found from the symmetry group of the network, as done in \cite{dellarossa:2020} for the orbital case and in \cite{siddique:2018} for the equitable single-layer case. For directed networks of two specific kinds (detailed in [9]), the matrix $T$ can be constructed as described in Sec. 2.3 in [9].
The key variational equation is reported here in compact form for ease of reference:
\begin{equation}\label{eq:variational}
	\dot\deta = [\rho_1 + \rho_2(B^k)] \deta .
\end{equation}
This equation describes the perturbation dynamics, by separating that along the synchronous manifold (described by the first $Q$ components $\eta_i$) from that transverse to it (described by the last components $\eta_i$, $i\in [Q+1, N]$).
We remark that the term $\rho_1$ in Eq. \eqref{eq:variational} is a diagonal matrix, which relates $\dot{\eta}_j$ only to $\eta_j$. By contrast, $\rho_2$ relates $\dot{\eta}_j$ also to the other perturbation components through the matrix $B^k$. Therefore, an inspection of the sub-blocks of $B^k$ allows to quickly check whether there is coupling between the dynamics of perturbations $\eta_i$ and $\eta_j$. 
To better illustrate this concept, let us consider the undirected, weighted network with $N=11$ nodes, $L=2$ kind of connections, and $Q=4$ clusters ($C_1,C_2,C_3,C_4$) shown in Fig. \ref{fig:netwTmatrix}, panel A, with nodes color coded to indicate the clusters they belong to. Note that the partition of the network nodes is equitable and not orbital \cite{siddique:2018}. Notice also the presence of a delay $\delta_2$ in the connection between nodes 5 and 6.

Panel C shows the structure of the matrices $T$ (left) and $B^1$ (right) for this network. Notice that matrix $B^2$ has the same structure as $B^1$, whose gray blocks contain only 0 entries. The upper-left $Q\times Q$ block is related to the perturbation dynamics along the synchronous manifold. Each white sub-block in the lower-right $(N-Q)\times (N-Q)$ sub-matrix $B^1_{N-Q}$ (with dashed black borders) describes the perturbation dynamics transverse to the synchronous manifold, thus is associated with loss of synchronization, either transient or permanent depending on the cluster stability. For instance, the $1 \times 1$ yellow (or blue or red) sub-block, is related to cluster $C_4$ (or $C_3$ or $C_1$, respectively), as pointed out in the corresponding row in matrix $T$, and describes the dynamics of the perturbation component $\eta_{11}$ (or $\eta_5$ or $\eta_{10}$, respectively); similarly, the $4 \times 4$ multi-color sub-block corresponds to clusters $C_1,C_2,C_3$. We remark that the structure of this sub-block implies that $\dot{\eta}_6, \dot{\eta}_7, \dot{\eta}_8, \dot{\eta}_9$ depend on $\eta_6, \eta_7, \eta_8, \eta_9$ but not on the other perturbations. Each transverse sub-block has an associated Maximum Laypunov Exponent (MLE) $\Lambda_i$, which can be studied independently from each other.

The stability of each cluster $C_q$ related to one or more sub-blocks depends on the maximum MLE $\Lambda_{C_q}$ among those associated to these sub-blocks: if $\Lambda_{C_q}$ is negative, the cluster $C_q$ is stable, otherwise it is unstable.
In the example, we computed the MLE associated to each sub-block: $\Lambda_1$ (blue sub-block), $\Lambda_2$ (multi-color sub-block), $\Lambda_3$ (red sub-block) and $\Lambda_4$ (yellow sub-block). The stability of $C_4$ depends on the sign of $\Lambda_{C_4} = \Lambda_4 = \max \{\lambda_{11}\}$ (i.e., the maximum component of the vector $\lambda_{11}$), whereas the stability of $C_1$ depends on the sign of $\Lambda_{C_1} = \max\{\Lambda_2, \Lambda_3\}$, the stability of $C_2$ depends on the sign of $\Lambda_{C_2} = \Lambda_2$ and the stability of $C_3$ depends on the sign of $\Lambda_{C_3} = \max\{\Lambda_1, \Lambda_2\}$.

Notice that the structure of the matrix $B^1$ allows us to state something more about the cluster stability. Indeed, the red cluster is related to two sub-blocks: the $1 \times 1$ red sub-block and the $4 \times 4$ multi-color sub-block. This means what follows: it is possible for the red cluster to undergo isolated desynchronization (see panel E) if the MLE $\Lambda_3$ becomes positive, while if the MLE $\Lambda_2$ becomes positive, red, blue, and green clusters become unstable together (see panel D).


This example clearly shows that the stability of each cluster in a subset of \textit{intertwined clusters} \cite{pecora:2014} may depend on the stability of the other clusters that belong to the same subset, but is decoupled from the clusters outside of the subset. Therefore, intertwined clusters can lose synchronization without causing a loss of synchronization in the clusters outside the subset, as for the yellow cluster in panel D.

\begin{figure}[b!]
	\centering
	\includegraphics[width=1\columnwidth]{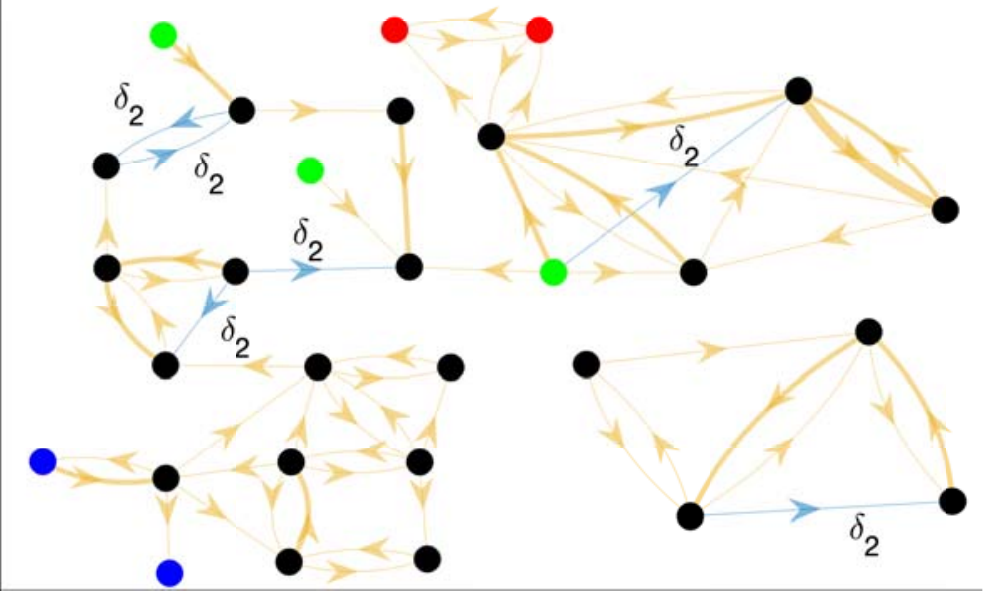}
	\caption{Macaque cortical connectivity network: $N = 29$ nodes, $M = 2$ node models, $L=2$ synapse models. Trivial clusters are black. Nodes of the same (non-black) color belong to the same cluster.}
	\label{fig:macaque_circuit}
\end{figure}

We apply the proposed method to a directed network (shown in Fig. \ref{fig:macaque_circuit}) composed of $N = 29$ nodes, each one representing one of the 91 areas of the macaque cerebral cortex \cite{markov:2013,markov:2014}. The neuron models that represent each area are of $M = 2$ kinds: 28 nodes are of kind $i=1$ and one node (corresponding to area V1) is of kind $i=2$, which is due to this one node receiving a visual input \cite{chaudhuri:2015}.
The nodes are connected through $L=2$ kinds of chemical excitatory synapses: one (for $k=1$) that transmits undelayed signals with $\delta_1 = 0$ (in yellow), one (for $k=2$) with delay $\delta_2 > 0$ (in blue).
The overall network is modeled by using the neuron and synapse equations described in Sec. 3 of [9] and the coupling matrices $A^1$ and $A^2$ provided in [9] (dataset S1). The measured connection weights \cite{markov:2014}, which range between 0 and 0.7636, have been quantized on four levels (0, 0.1, 0.5, 1) by replacing each original weight with  the closest one according to the Euclidean distance. After that, physical connections with length lower than 20 mm have been considered instantaneous (i.e., of kind $k=1$) and the corresponding quantized weights have been stored in the matrix $A^1$, whereas those longer than 20mm have been considered delayed (i.e., of kind $k=2$) and the corresponding quantized weights have been stored in the matrix $A^2$. These quantizations are justified by the fact that exact values for the coupling strengths and the delays reported in the literature are inevitably subject to measurement noise, and by the fact that, as we will see, they lead to the observation of functional mechanisms which are in agreement with physiological data, despite our simplifications.

The network non-trivial equitable clusters (consisting of more than one node) are listed in Tab. \ref{tab:macaque_clusters}. The same information is provided in Fig. \ref{fig:macaque_circuit}, where nodes of the same color (excluding black) belong to the same cluster: green for $C_1$, red for $C_2$ and blue for $C_3$. All nodes in trivial orbits are colored black. Obviously, the presence of a large number of trivial clusters does not mean that the corresponding areas are independent: they are densely connected, as evidenced in Fig. \ref{fig:macaque_circuit}, but they cannot be exactly synchronized.

\begin{table}[t!!]
	\centering
	\caption{ECs of the macaque cortical network}
	\label{tab:macaque_clusters}
	\begin{tabular}{cc}
		Cluster index & nodes in cluster \\
		\hline
		$C_1$ & 4;21;25 \\
		$C_2$ & 8,16 \\
		$C_3$ & 9,19 \\
		\hline
	\end{tabular}
\end{table}

Despite the rough quantizations applied to synaptic weights and delays, the clusters displayed in Table \ref{tab:macaque_clusters} are consistent with some previously reported physiological findings. For instance, cluster $C_2$ contains the nodes corresponding to visual areas 8l and 9/46v in the prefrontal cortex, which are known to be physically close and with similar connections \cite{markov:2014,goulas:2015}. The same holds for cluster $C_3$, which contains the nodes corresponding to the posterior and anterior portion of the inferotemporal cortex (TEO and TEpd, respectively).
\begin{figure}[!t]
	\centering
	\includegraphics[width=1\columnwidth]{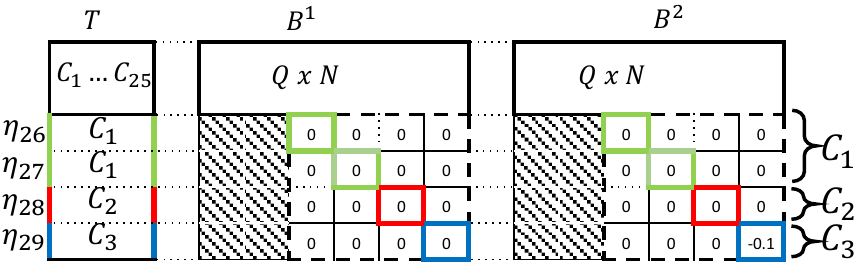}
	\caption{Structure of the matrices $T$, $B^1$, $B^2$, and $B^3$ for the macaque cerebral cortex. The gray blocks correspond to 0 entries.}
	\label{fig:B_mac}
\end{figure}

\begin{figure}[b!]
	\centering
	\includegraphics[width=1\columnwidth]{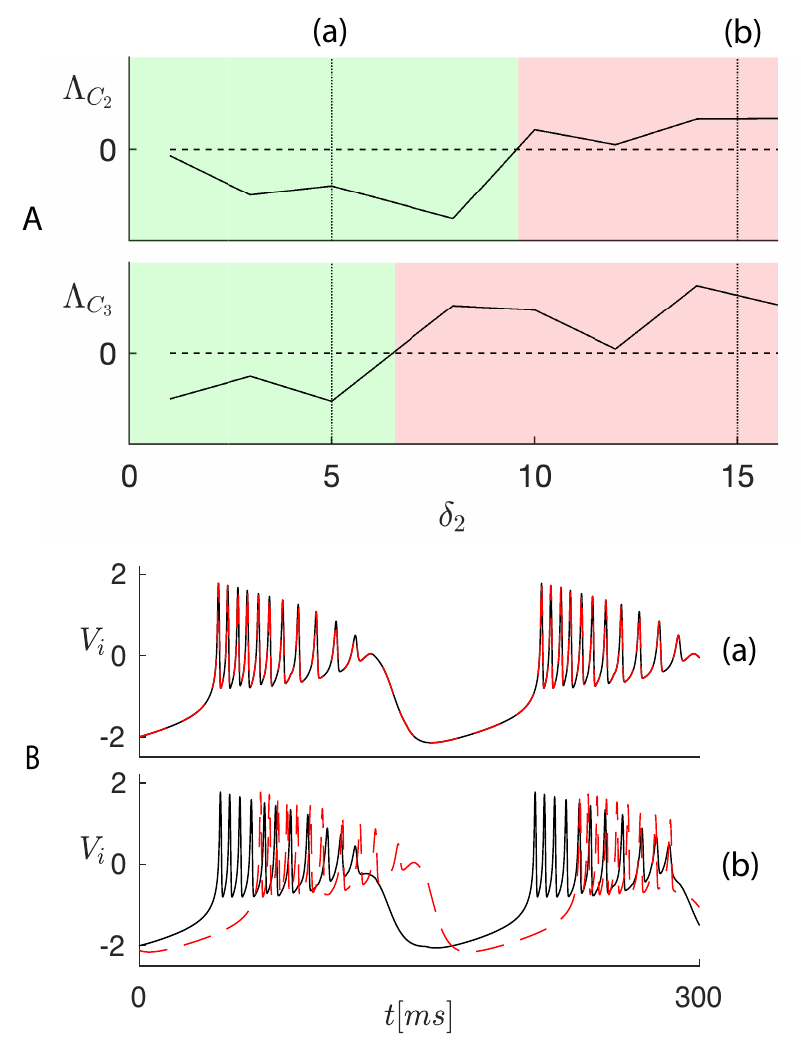}
	\caption{(A) MLE $\Lambda_{C_q}$ of each cluster $C_q$ ($q=2,3$) vs. coupling delay $\delta_2$, for the macaque cortical connectivity network. Horizontal dashed lines: edge of stability. Vertical dotted lines: $\delta_2$ values corresponding to the time plots in panel B. (B) Time plots $V_i(t)$ for different values of $\delta_2$ (5 ms (a), 15 ms (b)) for cluster $C_3$.}
	\label{fig:macaque}
\end{figure}

The directed connections originate from or go to trivial clusters only, therefore the cluster stability could be analyzed through the proposed approach.
Fig. \ref{fig:B_mac} shows the structure of the matrices $T$ (left) and $B^k$ (right) for this network.

In the block upper-triangular matrices $B^k$, the first $Q$ rows are related to the perturbation dynamics along the synchronous manifold. Each white sub-block in the lower-right $(N-Q)\times (N-Q)$ sub-matrix $B^k_{N-Q}$ describes the perturbation dynamics transverse to the synchronous manifold. 

If we analyze the matrices $B^k$ (related to the $k$-th connection type), we can see that: $B^1_{N-Q}$ (related to undelayed chemical excitatory synapses) has only zero entries, which implies that for the network with only these synapses the dynamics of each perturbation component $\eta_{k}$ depends only on $\eta_{k}$ through the term $\rho_1$ in Eq. \eqref{eq:variational};
$B^2_{N-Q}$ (related to delayed chemical excitatory synapses) has one $1 \times 1$ sub-block (with blue borders) related to cluster $C_3$, which means that for the network with only the delayed chemical excitatory synapses the dynamics of the perturbation component $\eta_{29}$ depends on $\eta_{29}$ through both $\rho_1$ and $\rho_2$ in Eq. \eqref{eq:variational}.
In summary, if we consider the whole network, with all kinds of synapses, the three clusters $C_1, C_2, C_3$ turn out to be not intertwined.

The stability analysis has been carried out by varying the delay $\delta_2$ between 0 and 16 ms (8 evenly spaced values). The neurons belonging to cluster $C_1$ do not receive any synaptic inputs, therefore the cluster transverse MLE is $\Lambda_{C_1} = 0$ for any value of $\delta_2$.
Figure \ref{fig:macaque}, panel A, shows the MLEs $\Lambda_{C_q}$ of the other clusters $C_q$ ($q=2,3$) versus the delay $\delta_2$. The green (red) regions in each plot $\Lambda_{C_q}(\delta_2)$ denote stability (instability) of the corresponding cluster $C_q$.

The vertical dotted lines mark the $\delta_2$ values corresponding to the time plots shown in panel B: $\delta_2 = 5$ ms (a) and $\delta_2 = 15$ ms (b). These plots display the first state variable $V_i$ of the neurons in cluster $C_3$. 
The panels show a window of 300 ms after a transient of 19.5 s. The breaking of this cluster is caused by a supercritical pitchfork bifurcation of cycles at each transition between the red and green regions, which generates two smaller stable trivial sub-clusters, each one producing one of the membrane voltages (black or red) in panel B, plot (b).

\begin{figure}[t!!]
	\centering
	\includegraphics[width=0.95\columnwidth]{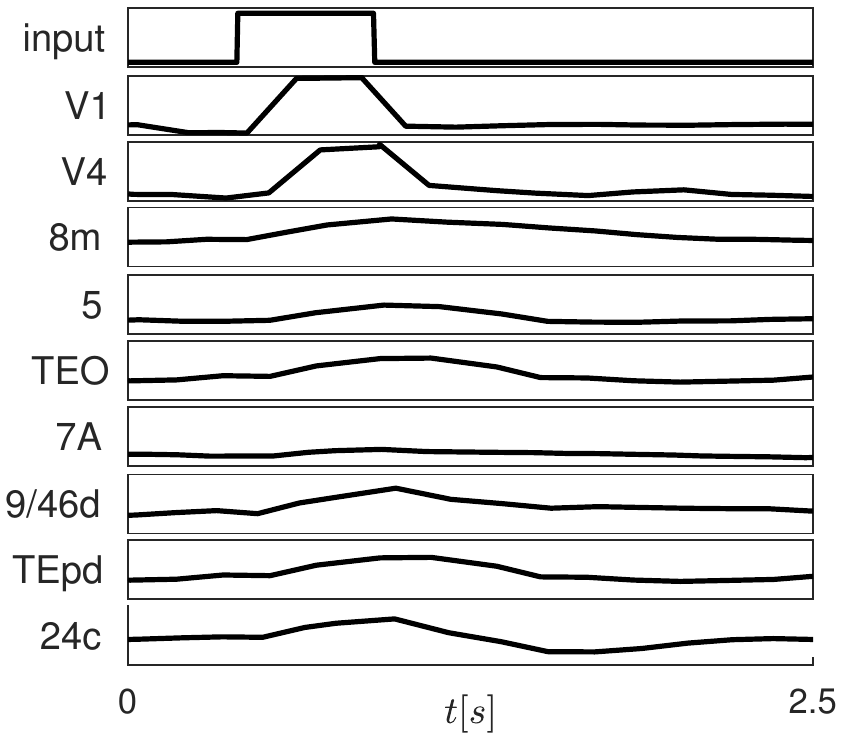}
	\caption{Time responses (firing rates) to a pulse-shaped input to area V1.}
	\label{fig:mac_time_resp}
\end{figure}

From Fig. \ref{fig:macaque} it clearly emerges that the two neurons in cluster $C_3$ display a phase lag for $\delta_2 = 15$ ms.
The synchronization of macaque visual cortex areas in response to visual stimuli has been observed in many experiments \cite{bosman:2012,chaudhuri:2015}. In particular, the areas 8l and 9/46v respond in a very similar way to visual inputs to area V1 \cite{chaudhuri:2015}. We thus set $\delta_2=5$ms in order to ensure synchronization of these two areas.

We proceeded to validate our model against the quantizations applied to the synaptic weights and axon delays, described before. To this end, following \cite{chaudhuri:2015} we simulated its response to a pulsed input to the primary visual cortex (area V1).
The response is propagated up the visual hierarchy, progressively slowing as it proceeds, as shown in Fig. \ref{fig:mac_time_resp}. Early visual areas, such as V1 and V4, exhibit fast responses. By contrast, prefrontal areas, such as 8m and 24c, exhibit slower decays to the standard firing rate, with traces of the stimulus persisting several seconds after stimulation. This is in agreement with the results shown in \cite{chaudhuri:2015} (compare with Fig. 3 in \cite{chaudhuri:2015}), which unveil a circuit mechanism for hierarchical processing of visual stimuli in the macaque cortex. Moreover, Fig. \ref{fig:mac_time_resp} evidences CS of the areas TEO and TEpd, corresponding to cluster $C_3$, as predicted by Fig. \ref{fig:macaque}.

The framework proposed in this paper is a fundamental step towards a method that fills the gap between the analysis of CS and networks of neurons.
The proposed method, based on a multi-layer network, allowed us to analyze the CS in the macaque cerebral cortex. A second example, which illustrates symultaneous loss of stability of intertwined clusters in a smaller network is provided in [9], showing an excellent agreement with biological measurements.

The results of this paper can be extended to study synchronization in any network characterized by different nodes, connections, and communication delays. As a final remark, we point out that the proposed model is completely deterministic and assumes that a reliable model of the network is available. These are quite strong modeling assumptions, since in real neuron networks the presence of noise is unavoidable and not always neuron and synapse models can be determined accurately. Despite this and despite the absence of information about the basins of attraction of stable clusters, our approach can provide useful information. For instance, in a real network CS will be approximate \cite{sorrentino:2016b}, not exact, as measured by high correlation values between the membrane potentials of the neurons/nodes belonging to a given stable cluster.
In this perspective, the patterns found with the proposed method are approximations to some more realistic solutions, which are characterized by higher complexity. Nature is quite far from determinism, therefore our analysis method is far from providing a general description of the dynamics of real neuron networks. This notwithstanding, it provides basic understanding of CS mechanisms, whose robustness can be checked by resorting to other less deterministic approaches.

\vspace{0.5cm}
The authors would like to express their sincere appreciation to Maurizio Mattia, Mauro Parodi and Lou Pecora for many useful inputs and valuable comments.


%

\end{document}